\documentclass[12pt]{amsart}
\usepackage{amsthm}
\usepackage{amsmath}
\usepackage{amssymb}

\newtheorem*{MDCT}{Morse Dynamical Characterization Theorem}
\newtheorem*{MCCT}{Morse Combinatorial Characterization Theorem}
\newtheorem*{TC}{Toeplitz Corollary}
\newtheorem*{MC}{Morse Corollary}
\newtheorem*{TDCT}{Toeplitz Dynamical Characterization Theorem}
\newtheorem*{TCCT}{Toeplitz Combinatorial Characterization Theorem}
\newtheorem*{theorem}{Theorem}

\def \sms {substitution minimal set}
\def \taubar {\overline\tau}
\def \mubar {\overline\mu}
\def \thetabar {\overline\theta}
\def \oneone {one-to-one}
\def \sub  {substitution}
\def \topconj {topologically conjugate}
\def \nlni {\newline\noindent}
\def \st   {such that }
\def \concat {concatenation}
\def  \s {\sigma}
\def \TMS {Toeplitz minimal set}
\def \MMS {Morse minimal set}
\def \ZZ {\mathbb Z}

\DeclareMathOperator{\cl}{cl}

\title[The Morse Minimal Set]{A Characterization of the Morse Minimal 
Set up to Topological Conjugacy}

\author{Ethan M. Coven}
\address{Ethan M. Coven,
Department of Mathematics,
Wesleyan University,
Middletown, CT 06459}
\email{ecoven@wesleyan.edu}
\author{Michael Keane}
\address{Michael Keane,
Department of Mathematics,
Wesleyan University,
Middletown, CT 06459}
\email{mkeane@wesleyan.edu}
\author{Michelle LeMasurier}
\address{Michelle LeMasurier,
Department of Mathematics,
Hamilton College,
Clinton, NY 13323}
\email{mlemasur@hamilton.edu}

\thanks{
Part of this work was done while the first author was Research
Professor at the Mathematical Sciences Research Institute (MSRI).
He thanks it for its hospitality.}

\subjclass[2000]{Primary 37B10}

\keywords{Morse minimal set, Toeplitz minimal set, topologically
conjugate}

\date{May 20, 2007}

\begin{document}

\begin{abstract}
We establish necessary and sufficient conditions for a dynamical
system to be \topconj\ to the \MMS, the shift orbit closure of the
Morse sequence, and conditions for topological conjugacy to the
closely related \TMS.
\end{abstract}

\maketitle

\section*{Introduction}

The celebrated {\em Prouhet-Thue-Morse sequence\/} is

  $$ 0110\ 1001\ 1001\ 0110 \dots $$

\noindent For an extensive bibliography, see \cite{AS99}. As there
may still be earlier occurrences of this sequence, it is not clear
what the name should be; for simplicity we adopt the common usage
{\em Morse sequence} in this article. For technical reasons we use
the two-sided Morse sequence

$$  m = \dots 10010110.01101001 \dots $$

\noindent Here the decimal point separates the negative indices from
the non-negative ones, and the left side is the right side reflected
across the decimal point.

The shift orbit closure

$$ M := \cl \{\s^n(m) : n \in \ZZ\}$$

\noindent of the two-sided Morse sequence in the space of all
bilaterally infinite zero-one sequences is called the {\em \MMS}.
The members of this set all have the same collection of {\em
blocks}\ , i.e. finite words. (In the terminology of formal language
theory, they have the same language.) It is well-known that the
members of the \MMS\ are characterized by the ``no $BBb$" property:
they contain no block of the form $BBb$, where $b$ is the first
letter of~$B$.

In this article, we characterize the \MMS \ {\em dynamically}, i.e.\
up to topological conjugacy. Therefore it will be useful to think of
the \MMS\ as a minimal symbolic dynamical system~$(M,\s)$. Along the
way we characterize, both {\em combinatorially}, i.e. up to
equality, and dynamically the closely related {\em \TMS},
denoted~$(T,\s)$. The name originates in \cite{JK69}, although there
are earlier appearences (see below); there is no name ambiguity here
as it was chosen to be descriptive in this reference.

Our analysis of the Morse and \TMS s is based on the fact that both
are symbolic minimal sets generated by substitutions of constant
length.  A {\em substitution\/} of constant length $r \ge 2$ is a
mapping
$$\theta : A \to A^r,$$
where $A$ is a finite set (the {\em alphabet}) with at least two letters.
Such a \sub\ determines a topological semi-conjugacy
$$ \thetabar : (A^{\ZZ}, \s) \xrightarrow{\text{into}}
(A^{\ZZ}, \s^r)$$
as follows.  For every $x = (x_n) \in A^\ZZ$ and for every $n \in \ZZ$,
$$[\thetabar(x)]_{nr,nr+1,\dots,nr+r-1} := \theta(x_n).$$
Then ${\overline{\theta^k}} =\thetabar^k$ for every $k \ge 1$. A
{\em \sms\ generated by $\theta$} is a minimal shift orbit closure
of a $\thetabar$-periodic point. For a general reference on
substitution minimal sets, see \cite{G}.

\smallskip

The Morse \sub\ is
$$\mu(0) := 01, \quad \mu(1) := 10$$
and the Toeplitz \sub\ is
$$\tau(0) := 01, \quad \tau(1) := 00.$$
Since both $\mu$ and $\tau$ are \oneone, so are $\mubar$ and
$\taubar$, and hence both are topological conjugacies from their
domains onto their images.

For the Morse \sub, $\mubar$ has four periodic points, all with
least $\mubar$-period~2, denoted $m_{a.b}$  ($a,b = 0,1$). $m_{a.b}$
is the unique $\mubar^2$-fixed point $x$ \st $x_{-1} = a$ and $x_0 =
b$. The \MMS\ is the shift orbit closure of any of these four
points.
The two-sided Morse sequence $m$ defined above is ~$m_{0.0}$. For
the Toeplitz \sub, $\taubar$ has two periodic points, $t_{0.0}$ and
$t_{1.0}$, with the same orbit closure. Both have least  $\taubar$-period~2.

\section*{Combinatorial Characterizations}

In this section we give combinatorial characterizations of the Morse
and Toeplitz minimal sets. The relevant property that relates them
is that the \TMS\ is an exactly two-to-one image of the \MMS\ via
the topological semi-conjugacy $\Phi$, given by the local rule
$$\varphi(u,v) := u+v+1 \mod 2.$$
In \cite{GH55} this is attributed to J.~C.~Oxtoby (p.~113).

The following result for the \MMS\ is well-known.

\begin{MCCT}
A bilaterally infinite sequence with entries from $\{0,1\}$ is in the
\MMS\ if and only if it contains no block of the form $BBb$, where
$b$ is the first letter of~$B$.
\end{MCCT}

\smallskip
\noindent {\bf Remark.}  The  ``no $BBb$" property is  what A. Thue
used to define the Morse sequence (\cite{Th06},\cite{Th12}, see also
\cite{Th77}). Thue's work  was unknown in the world of dynamics
until~1967, when G. A. Hedlund  learned  it from R. B\"{u}chi (see
\cite{He67}). That every member of the \MMS\ has  the no $BBb$
property was proved by M. Morse and Hedlund \cite{MH44} and that
these are the only points having it by W. H. Gottschalk and Hedlund
\cite{GH64}. Morse had proved some years earlier  that the Morse
sequence has the ``no $BBb$" property \cite{Mo38}; it is obvious to
us today that the Morse result implies the Morse-Hedlund result.

\begin{TCCT}
A bilaterally infinite sequence with entries from $\{0,1\}$
is in the \TMS\ if and only if it contains no block of the form $BB$
where the number of zeros in
~$B$
is even.
\end{TCCT}

\begin{proof}
Let $\Phi$ be the topological semi-conjugacy defined above and let
$\varphi$ be its local rule.  Since $\Phi(M) = T$, it suffices to
show that $\Phi(M)$ is the set of bilaterally infinite sequences with
entries from $\{0,1\}$ that contain no block of the form ~$BB$
where the number of zeros in
~$B$ is even.

Setting $\bar 0 = 1$ and $\bar 1 = 0$, $\Phi$ identifies  $x$
and~$x'$ if and only if  $x' = \bar x$, and  $\varphi$ identifies
$B$ and ~$B'$ if and only if  $B' = \bar B$. Therefore for every
block $B$ with entries from$\{0,1\}$, there is  a block ~$C$ of the
same length as~$B$ and with first letter~$c$, and
$$\varphi^{-1}(BB) = \{CCc,\bar{C}\bar{C}\bar{c}\} \text{ or }
  \{C\bar{C}c,\bar{C}C\bar{c}\}$$
according as the number of zeros in $B$ is even or odd.
\end{proof}

\section*{Dynamical Characterizations}In this section we give dynamical
characterizations of the Morse and Toeplitz minimal sets. Although
the results are stated for symbolic minimal sets, recall that if $X$
is a zero-dimensional, compact, metrizable space and if $S$ is a
homeomorphism of $X$, then $(X,S)$ is topologically conjugate to a symbolic
minimal set if and only if $S$ is minimal and expansive. This is our
justification for using the words ``dynamical characterizations."

For our purposes, the key property of the \MMS\ $(M,\s)$  is that
for every $k \ge 0$, every point in~$M$ can be written as a \concat\
of the $2^k$-blocks $\mu^k(0)$ and $\mu^k(1)$,
$$ \dots \mu^k(m_{-1})\mu^k(m_0)\mu^k(m_1)\dots$$
for some $m = (m_i) \in M$.
An analogous statement holds for the \TMS.

The dynamical characterizations of the Morse and \TMS s are
conceptually the same, but the details of the proof in the Toeplitz
case are simpler.  For that reason, we state and prove it first.

\begin{TDCT}
Let $(Y,\s)$ be a symbolic minimal set. Then $(Y,\s)$ is \topconj\
to the \TMS\  if and only if there exist $k\ge0$ and $2^k$-blocks
$C_0 \ne C_1$ \st\ some  (and hence every) point in~$Y$ can be
written as a \concat\ of~$C_0$ and ~$C_1$ that, thought of as a
bilaterally infinite sequence on letters~$C_0$ and ~$C_1$,  contains no
block of the form $BB$ where the number of $C_0$ in ~$B$ is even.
\end{TDCT}

\begin{proof}
(1) Suppose that $(Y,\s)$ is a symbolic minimal set that is
\topconj\ to the \TMS\ $(T,\s).$ There is a topological conjugacy
$\Psi:(T,\s) \to (Y,\s)$ that   is given by a local rule $\psi$ with
no memory and anticipation $2^k-1$ for some $k\ge1$. For every $t
\in T$ and for every $i \in \ZZ$,
$$[\Psi(t)]_i = \psi(t_i,t_{i+1},\dots,t_{i+2^k-1}).$$
Thus $\psi$ maps the $2^k$-blocks of $T$ to 1-blocks of $Y$, and for
every $j\ge1$, $\psi$ maps the  $(2^k+j-1)$-blocks of $T$ to
$j$-blocks of $Y$. Taking $j = 2^k$, $\psi$ maps the
$(2^{k+1}-1)$-blocks of $T$ to ${2^k}$-blocks of $Y$.

Let $y \in Y$.  Set $t := \Psi^{-1}(y)$ and write
$$t = \dots \tau^k(s_{-1})\tau^k(s_0)\tau^k(s_1)\dots$$
for some $s = (s_i) \in T$. An easy induction shows that for every
$\ell \ge 1$, $\tau^{\ell}(0)$ and $\tau^{\ell}(1)$ agree everywhere
except in their last places, so the initial $(2^{k+1}-1)$-blocks of
$\psi(\tau^k(0)\tau^k(0))$ and of $\psi(\tau^k(0)\tau^k(1))$ are
equal.  Note that since 11 does not appear in ~$T$,
$\psi(\tau^k(1)\tau^k(1))$ does not appear in~$y$.

So define $C_0$  to be the common initial $2^k$-block of
$\psi(\tau^k(0)\tau^k(0))$ and of  $\psi(\tau^k(0)\tau^k(1))$, and
$C_1$ to be the initial $2^k$-block of $\psi(\tau^k(1)\tau^k(0))$.
If $C_0 = C_1$, then  $\Psi(T)$  would be finite.
Thus $C_0 \ne C_1$ and $y$ can be written as a \concat\ of $C_0$ and
$C_1$. Thinking of the \concat\ as a bilaterally infinite sequence with
letters $C_0$ and~$C_1$,  it contains no block of the form $BB$
where the number of~$C_0$ in ~$B$ is even, for otherwise $s$ would
contain a block of the form $BB$ where the number of zeros in ~$B$
is even.

\smallskip

(2) Conversely, assume that there are $2^k$-blocks $C_0 \ne C_1$ \st
every point in~$Y$ can be written as a \concat\ of~$C_0$ and ~$C_1$
that, thought of as a bilaterally infinite sequence on letters~$C_0$ and
~$C_1$,  contains no block of the form $BB$ where the number
of~$C_0$ in ~$B$ is even.

Let $Y'$ be the set of points in~$Y$ \st\ for every $n \in \ZZ$,
$$y_{2^kn}y_{2^kn+1}\dots
y_{{2^kn}+2^k-1} = C_0 \text{ or } C_1,$$
and let $T'$ be the corresponding subset of $T$,
i.e. the set of points $t \in T$ \st\
for every $n \in \ZZ$,
$$t_{2^kn}t_{2^kn+1}\dots
t_{{2^kn}+2^k-1} = \tau^k(0) \text{ or } \tau^k(1).$$
Then the map $\Upsilon : Y' \to T'$, defined by
$$ C_0 \mapsto \tau^k(0) \qquad C_1 \mapsto \tau^k(1),$$
is a topological conjugacy of $(Y',\s^{2^k})$
onto $(T',\s^{2^k})$.
Extend $\Upsilon$ to all of $Y$ by defining
$$\Upsilon := \s^j \circ \Upsilon \circ \s^{-j} \text{ on }
\s^j(Y'), \quad j=1,2,\dots,2^k-1.$$ $\Upsilon$ is well-defined and
\oneone, and hence a topological conjugacy of $(Y,\s)$ onto
$(T,\s)$, provided that $ \{Y',\s(Y'),\dots,\s^{2^k-1}(Y')\}$ is
pairwise disjoint.

If not, there exists  $\ell$,  $1 \le \ell \le 2^k-1$,  \st\ $Y'
\cap \s^\ell(Y') \ne \varnothing$. Since every point of $Y'$ can be
written as a \concat\ of~$C_0$ and ~$C_1$, and  $C_1$  never follows
itself, $C_0$ must overlap both $C_0$ and $C_1$,  and $C_1$ must
overlap $C_0$, all overlaps starting at the same place.  It is then
straightforward to show that $C_0 = C_1$. Since $C_0 \ne C_1$, $
\{Y',\s(Y'),\dots,\s^{2^k-1}(Y')\}$ must be  pairwise disjoint and
the proof is complete.
\end{proof}

\noindent {\bf Remark.} The number of symbols in the symbolic
minimal set can be any finite number. However, we will prove in the
next section that a \sms\ that is \topconj\ to the \TMS\ must be on
an alphabet of at most three symbols.

\smallskip

  The proofs of the Toeplitz and Morse dynamical theorems   are
conceptually the same, differing only in the details.  We will
therefore  only sketch the proof of the latter, indicating the
differences.

\begin{MDCT}
Let $(X,\s)$ be a symbolic minimal set. Then $(X,\s)$ is \topconj\
to the \MMS\  if and only if there exist $k\ge0$ and $2^k$-blocks
$C_0 \ne C_1,C_0'$, and~$C_1'$  \st some (and hence every) point
in~$X$ can be written as a \concat\ of~$C_0, C_1, C_0'$, and~$C_1'$
satisfying \nlni (i) every second block in the \concat\ is a
sequence of~$C_0$ and~$C_1$ that, thought of as a bilaterally infinite
sequence with letters $C_0$ and~$C_1$, contains no block the form
$BBb$, where $b$ is the first letter of~$B$. \nlni (ii) the gaps are
filled in according to the following nearest neighbor rule.
$$ C_0  \text{ goes between } C_1 \text{ and } C_1$$
$$ C_1  \text{ goes between } C_0 \text{ and } C_0$$
$$C_0' \text{ goes between } C_1 \text{ and } C_0$$
$$C_1' \text{ goes between } C_0 \text{ and } C_1.$$

\end{MDCT}

\begin{proof}
(1) Suppose that $(X,\s)$ is a symbolic minimal set that is
\topconj\ to the \MMS\ $(M,\s).$ There is a topological conjugacy
$\Phi:(M,\s) \to (X,\s)$ that  is given by a local rule $\varphi$
with no memory and anticipation $2^k$ for some $k\ge0$. $\varphi$
maps the  $2^{k+1}$-blocks of $M$ to ${2^k}$-blocks of $X$.

Every point of $M$ can be written as a \concat\ of the $2^k$-blocks
$\mu^k(0)$ and $\mu^k(1)$ and of the $2^{k+1}$-blocks
$\mu^{k+1}(0) = \mu^k(0)\mu^k(1)$ and $\mu^{k+1}(1) = \mu^k(1)\mu^k(0)$.

Define
$$C_0 := \varphi(\mu^k(0)\mu^k(1))$$
$$C_1 := \varphi(\mu^k(1)\mu^k(0))$$
$$C_0' := \varphi(\mu^k(0)\mu^k(0))$$
$$C_1' := \varphi(\mu^k(1)\mu^k(1)).$$
Look at the set  $M_0$  of points in~ $M$  such that the
$2^{k+1}$-blocks starting at multiples of ~$2^{k+1}$  are
$\mu^{k+1}(0)$  and  $\mu^{k+1}(1)$.  $\Phi$  is a topological
conjugacy of  $(M_0,\sigma^{2^{k+1}})$  onto
$(\Phi(M_0),\sigma^{2^{k+1}})$.  Then  $C_0 \ne C_1$, for otherwise
$\Phi(M_0)$ and hence  $\Phi(M)$  would be finite. That
$C_0,C_1,C_0'$, and $C_1'$ satisfy conditions (i) and~(ii) follows
from the definitions.

\smallskip

(2) Conversely, assume that there are $2^k$-blocks $C_0\ne C_1,C_0'$,
and~$C_1'$ satisfying conditions (i) and~ (ii).

Let $X'$ be the set of points in~$X$ \st for every $n \in \ZZ$,
$$x_{2^kn}x_{2^kn+1}\dots
x_{{2^kn}+2^k-1} = C_0, C_1,C_0',
\text{ or } C_1'$$ and the \concat\ satisfies
conditions (i) and~ (ii);
and let $M'$ be the corresponding subset of $M$.

Since $C_0$ can be followed only by $C_1$ or $C_1'$, and $C_1$ can
be followed only by $C_0$ or $C_0'$, the map $\Xi$, defined by
$$ C_0,C_0' \mapsto \mu^k(0) \qquad C_1,C_1' \mapsto \mu^k(1),$$
is a well-defined topological conjugacy of $(X',\s^{2^k})$
onto $(M',\s^{2^k})$.
Extend $\Xi$ to all of $X$ by defining
$$\Xi := \s^j \circ \Xi \circ \s^{-j} \text{ on }
\s^j(X'), \quad j = 1,2,\dots,2^k-1.$$

\noindent $\Xi$ is well-defined and \oneone, and hence a topological
conjugacy of $(X,\s)$ onto $(M,\s)$, provided that $
\{X',\s(X'),\dots,\s^{2^k-1}(X')\}$ is pairwise disjoint. If not,
then as in the proof of the Toeplitz Theorem, $C_0 = C_1$.
\end{proof}

\noindent {\bf Remark.} One might guess from the proof that
$C_0=C_0'$ or $C_1=C_1'$ is possible. Taking $(X,\s)$ to be the
\MMS\ itself and $k=0$ shows this guess to be correct.

\medskip

As corollaries of the dynamical characterizations, we show that,
subject to some harmless restrictions on the \sub, a \sub\ of
constant length that generates a \sms\ \topconj\ to the Morse (resp.
Toeplitz) minimal set must be defined on an alphabet of at most six
(resp. three) letters. To state and prove these corollaries we use
the following terminology from the theory of directed graphs.

For a \sub\ $\theta$ of constant length defined on alphabet ~$A$,
let ~$G(\theta)$ be the directed graph  with vertices~$A$ and an arc
from~$a$ to~$b$ if $b$ appears in ~$\theta(a)$.   A directed graph is
{\em strongly connected\/} if for every ordered pair $(a,b)$ of
vertices, there is a directed path from~$a$ to ~$b$.  A strongly
connected directed graph with vertices~$A$ is either {\em primitive}
(there exists $K \ge 1$ such that for every ordered pair $(a,b)$ of
vertices,  there is a directed path of length~$K$ from~$a$ to ~$b$)
or it has a {\em period\/} $\ell \ge 2$ ($A$ can be written as a
disjoint union $A_0 \cup A_1 \cup \dots \cup A_{\ell -1}$ \st there
is an arc from a vertex in~$A_i$ to one in~$A_j$ only if $j \equiv
i+1 \mod \ell$).

\begin{TC}
Let $(Y,\s)$ be a \sms\ generated by a \oneone, primitive \sub\ ~
$\theta$ of constant length.  If  $(Y,\s)$ is \topconj\ to the \TMS,
then \nlni (a) the length of $\theta$ is a power of~2, \nlni (b) the
alphabet $A$ of $\theta$ has at most three letters.
\end{TC}

\begin{proof}
(a) An infinite \sms\ generated by a \sub\ of constant length~$r$ is
topologically semi-conjugate to the ``+1" map on the $r$-adic
integers.  Thus the set of prime factors of the length of the \sub\
is a topological conjugacy invariant.

\noindent (b)  By the Toeplitz dynamical theorem,
there exist  $k\ge0$ and $2^k$-blocks $C_0 \ne C_1$ such that every
point in ~$Y$ can be written as a \concat\ of~$C_0$ and~$C_1$ that
belongs to the \TMS\ on alphabet~$\{C_0,C_1\}$.  By replacing  $C_0$
and~$C_1$ by longer blocks or by replacing $\theta$ by a power, we
may assume that the length of~$\theta$ is~$2^k$.

Let $y \in Y$.  Since
$$Y = \thetabar(Y) \cup \s(\thetabar(Y)) \cup \dots \cup 
\s^{2^k-1}(\thetabar(Y)),$$
$y \in \s^j(\thetabar(Y))$ for some $j$, $0 \le j \le 2^k-1$.
Therefore $y$ can be written as a \concat\ of the $2^k$-blocks
$\theta(i)$, $i \in A$, all starting at places congruent to ~$j$
modulo~ ${2^k}.$ Since 00, 01, and 10  appear  in ~$T$ but 11 does
not,  $C_0C_0$, $C_0C_1$, and $C_1 C_0$ appear in~$Y$ but $C_1C_1$
does not. Thus every appearance of every $\theta(i)$ in the \concat\
is as a subblock of $C_0C_0$, $C_0C_1$, or $C_1 C_0$, starting at
place~$j$ of the left-hand $2^k$-block. So there are at most three
distinct $2^k$-blocks~$\theta(i)$.  Since $\theta$ is \oneone, $A$
has at most three letters.
\end{proof}

For example, the primitive, \oneone\ \sub
$$0\mapsto 12, \quad 1\mapsto 02, \quad 2\mapsto 10$$
generates a minimal set that
is \topconj\ to the \TMS\ via $C_0 = 21, C_1 = 00$.

\begin{MC}
Let $\theta$ be a \oneone, primitive \sub\ of constant length that
generates a unique minimal set  $(X,\s)$.  If  $(X,\s)$ is \topconj\
to the \MMS, then \nlni (a) the length of $\theta$ is a power of~2,
\nlni (b) the alphabet $A$ of $\theta$ has at most six letters.
\end{MC}

\begin{proof}
In the proof of the Toeplitz Corollary, replace $C_0$ and $C_1$  by
$C_0, C_1, C_0'$, and $C_1'$ from the Morse dynamical
theorem.  Of the sixteen ordered 2-tuples formed by
them, only $C_0C_1,C_0C_1',C_1 C_0,C_1 C_0',C_0'C_0,$ and $C_1'C_1$
can appear in ~$X$.
\end{proof}

{\bf Remark.} It is straightforward to show that an infinite \sms\
generated by a substitution of constant length can also be generated
by a primitive  substitution of  constant length, although the
primitive \sub\ may be of longer length.  If the  primitive \sub\ is
not \oneone, then identifying symbols with the same image yields a
primitive \sub\ on fewer symbols, and the \sms\ generated by it  is
\topconj\ to the original.  Thus the assumptions on the substitution
in the corollaries are harmless.

\section*{Dynamical systems conjugate to \sms s}

In this section we characterize those dynamical systems that are
topologically conjugate to \sms s generated by \sub s of constant
length.  As noted earlier, ``symbolic minimal'' is a topological
conjugacy invariant.

\begin{theorem}
Let $(X,\s)$ be an infinite  symbolic minimal set and let $r\ge2$.
Then $(X,\s)$ is \topconj\ to a \sms\   generated by a primitive
substitution of a constant length $r$ if and only if there is a
compact, $\s^r$-invariant, proper subset  $X'$ of $X$ such that
$(X,\s)$ is topologically conjugate to $(X',\s^r)$.
\end{theorem}

\begin{proof}

One implication (if $(X,S)$ is \topconj\ to a \sms\ ...)   is
\cite{BDM04}, Theorem 2.1, which is attributed there to B. Moss\'e.

So suppose that there is a compact, $\s^r$-invariant, proper subset
$X'$ of $X$ such that  $(X,\s)$ is topologically conjugate to
$(X',\s^r)$. Let $\Phi$ be a topological conjugacy of $(X,\s)$ onto
$(X',\s^r)$, given by a local rule $\varphi$, which we may assume
has no memory and anticipation  $m \ge 0$. Thus $\varphi$ maps the
$(m+1)$-blocks of $(X,\s)$ to $1$-blocks of $(X',\s^r)$.  By adding
superfluous variables, we may assume that  $\varphi$ maps the
$(m(r-1)+1)$-blocks of~$(X,\s)$ to $1$-blocks of~$(X',\s^r)$, i.e.
$r$-blocks of~$(X,\s)$. As before, $\varphi$ maps the
$(mr+1)$-blocks of~ $X$ to $(mr+r)$-blocks of~$X$.

Code the $(mr+1)$-blocks of ~$X$,  and let $(Y,\s)$ denote the
coded system with  $(X',\s^r)$ coding to~$(Y',\s^r)$. The
topological conjugacy  $\Phi$ codes to a topological conjugacy
~$\Theta$ with no memory and no anticipation.

Let  $A$  be the alphabet of~$(Y,\s)$.
   The local rule $\theta$ of ~$\Theta$
is a \sub\ of constant length~$r$. We show that $\theta$ is
primitive and hence generates a unique  minimal set, which must be
topologically conjugate to~$(X,\s)$.

If $G(\theta)$ is not strongly connected, then not all the symbols
in  ~$A$ appear in~$Y$. And if $G(\theta)$ has a period $\ell \ge
2$, then again not all the symbols from~$A$  appear in~$Y$.
Similarly, a primitive \sub\ generates a unique \sms.
\end{proof}

\end{document}